\documentclass[12pt]{article}

\usepackage{amsfonts}
\usepackage{amssymb}
\usepackage{amsthm}
\usepackage{amsmath}
\usepackage{graphicx}
\usepackage{color}

\usepackage{amscd}
\usepackage{multirow}

\newcommand{\CC}{{\mathbb C}}
\newcommand{\C}{{\mathbb C}}

\newcommand{\RR}{{\mathbb R}}

\newcommand{\ZZ}{{\mathbb Z}}
\newcommand{\Z}{{\mathbb Z}}
\newcommand{\NN}{{\mathbb N}}

\newcommand{\vol}{{\operatorname{vol}}}
\newcommand{\rank}{{\operatorname{rank}}}

\newtheorem{theorem}{Theorem}[section]

\newtheorem{example}[theorem]{Example}
\newtheorem{lemma}[theorem]{Lemma}
\newtheorem{cor}[theorem]{Corollary}
\newtheorem{remark}[theorem]{Remark}

\title{Exponential growth of rank jumps for $A$--hypergeometric systems.}

\author{Mar\'ia-Cruz Fern\'andez-Fern\'andez \thanks{Partially supported by MTM2010-19336 and FEDER, and Junta de Andaluc\'ia under grants FQM-5849, FQM333.
E.mail address: {\tt mcferfer@us.es}}\\
Departamento de \'{A}lgebra \\ Universidad de Sevilla}

\date{January 24, 2012}

\begin{document}

\maketitle

\begin{abstract}
The dimension of the space of holomorphic solutions at nonsingular points (also called the
holonomic rank) of a $A$--hypergeometric system $M_A (\beta )$ is
known to be bounded above by $ 2^{2d}\operatorname{vol}(A)$
\cite{SST}, where $d$ is the rank of the matrix $A$ and $\vol (A)$
is its normalized volume. This bound was thought to be very
vast because it is exponential on $d$. Indeed, 
all the examples we have found in the literature verify that
$\operatorname{rank}(M_A (\beta ))<2 \vol (A)$. We construct here, in a very elementary way, 
some families of matrices $A_{(d)}\in \ZZ^{d \times n}$ and
parameter vectors $\beta_{(d)} \in \C^d$, $d\geq 2$, such that
$\rank ( M_{A_{(d)}} (\beta_{(d)} ))\geq a^d \vol(A_{(d)})$ for
certain $a>1$.
\end{abstract}

\section{Introduction}
Let $A=(a_{ij})=(a_1 \; a_2 \cdots a_n)$ be a full rank matrix with columns $a_j \in \ZZ^d$ and $d\leq n$.
Following Gel'fand, Graev, Kapranov and Zelevinsky (see \cite{GGZ} and \cite{GKZ}) we can define
the $A$--hypergeometric system with parameter $\beta \in \C^d$ as the left ideal $H_A (\beta )$
of the Weyl algebra $D=\C [x_1 ,\ldots ,x_n ]\langle \partial_1 ,\ldots ,\partial_n \rangle$ generated by
the following set of differential operators:

\begin{equation}
\Box_u := (\prod_{i: u_i >0} \partial_{i}^{u_i})- (\prod_{i: u_i < 0}\partial_{i}^{-u_i}) \; \;\; \; \mbox{ for all }
u\in \ZZ^n \mbox{ such that } Au=0 \label{Toric-operators}
\end{equation} and
\begin{equation}
E_i - \beta_i := \sum_{j=1} ^n a_{ij}x_j\partial_j -\beta_i \; \; \;
\; \mbox{ for } i=1,\ldots, d \label{Euler-operators}
\end{equation}

\vspace{.3cm}

The operators given in (\ref{Toric-operators}) generate the
so-called toric ideal $I_A \subseteq
\CC[\partial_1,\ldots,\partial_n] $ associated with $A$ and the $d$ operators given in (\ref{Euler-operators}) are called the 
Euler operators associated with the pair $(A ,\beta)$. The hypergeometric $D$--module associated with the pair $(A,\beta )$
is the quotient $M_A(\beta)=D /D H_A (\beta )$. It is a holonomic
$D$--module for any pair $(A,\beta)$ as above (see \cite{GKZ},
\cite{Adolphson}). In particular, the space of holomorphic solutions
of $M_A (\beta)$ at a nonsingular point has finite dimension. This
dimension or, equivalently, the holonomic rank of $M_A (\beta )$
equals the normalized volume $\vol_{\ZZ A} (A)$ of the matrix $A$
(see (\ref{volume})) when either $I_A$ is Cohen-Macaulay or $\beta$
is generic (see \cite{GKZ}, \cite{Adolphson}, \cite{SST}).

The first example of a pair $(A,\beta)$ for which $\rank (M_A (\beta
))>\vol_{\ZZ A}(A)$ was described in \cite{ST98} (see Example
\ref{rank-jump-d2}). A complete description of the case $d=2$
appears in \cite{CDD99}, revealing that in this case the rank of
$M_A (\beta)$ can be only $\vol_{\ZZ A} (A)$ (the generic value) or
$\vol_{\ZZ A} (A) +1$ (the exceptional value).

In general it is known that $\rank (M_A (\beta ))\geq \vol_{\ZZ A}(A)$ for all $\beta$ \cite{SST, MMW}.
In fact, it is proved in \cite{MMW} that the map $\beta \in \CC^d \mapsto \rank (M_A (\beta))$ is upper
semi--continuous in the Zarisky topology and they also provide an
explicit description of the exceptional set
$$\varepsilon (A)=\{\beta \in \CC^d :\; \rank (M_A (\beta ))> \vol_{\ZZ A} (A) \}$$
that turns out to be an affine subspace arrangement with codimension at least $2$. Previous descriptions of the
exceptional set in particular cases appear in \cite{CDD99, Mat01, Saito, Mat03}.

If for a fixed matrix $A$ we have that $j_A (\beta )=\rank (M_A
(\beta ))- \vol_{\ZZ A} (A) >0$ then it said that the $A$--hyper\-geometric
system has a \emph{rank jump} of $j_A (\beta )$ at $\beta$ or that $\beta$ is a rank jumping parameter for
$A$.

The paper \cite{MW} provides the first family of hypergeometric
systems with rank jump greater than $2$. Indeed, they construct a
family of pairs $(A_{(d)} ,\beta_{(d)} )$ with $A_{(d)} \in \ZZ^{d\times 2d}$
and $\beta_{(d)} \in \CC^d$ such that $j_{A_{(d)}}(\beta_{(d)})= d-1$. However, for this family $\vol_{\ZZ A_{(d)}}(A_{(d)})=d+2$ and thus
$$\dfrac{\rank (M_{A_{(d)}} (\beta_{(d)} ))}{\vol_{\ZZ A_{(d)}}(A_{(d)})} =  2- \dfrac{3}{d+2}<2$$

More recently, in \cite{Berkesch} a general combinatorial formula is
provided for the rank jump $j_A (\beta )$ of the $A$--hypergeometric
system at a given $\beta$. However, the formula is very complicated
and, in fact, all the examples included in \cite{Berkesch} verify
that $\rank (M_A (\beta))< 2 \vol_{\ZZ A} (A)$ as well. Previous
computations of $j_A (\beta)$ in particular cases appear for example in 
\cite{CDD99}, \cite{Saito}, \cite{Oku}.

In the case when the toric ideal is standard homogeneous, the
following upper bound for the holonomic rank of a hypergeometric
system is proved in \cite{SST}:
$$
\operatorname{rank}(M_A (\beta )) \leq 2^{2d}\vol_{\ZZ
A}(A)\label{upper-bound}$$ However, it is mentioned in \cite[p. 159]{SST} that this
upper bound is most likely far from optimal and that it would be
desirable to know whether the ratio $\rank ( M_A (\beta))/\vol_{\ZZ A}
(A)$ can be bounded above by some polynomial function in $d$.
Here we provide a very elementary construction of some families of
hypergeometric systems for which the ratio $\rank ( M_{A} (\beta
))/\vol_{\ZZ A}(A)$ is exponential on $d$, giving a negative answer
to this last question.

Moreover, for one of the families constructed the dimension of Laurent polynomial solutions is lower 
than the rank jump (see Remark \ref{remark-polynomial-solutions}). This is in contrast with the general 
observation in the examples found in the literature (see for example \cite{MW}).

I am grateful to Christine Berkesch for many helpful conversations about her paper \cite{Berkesch}.

\section{Construction of the examples.}

Recall that the normalized volume of a full rank matrix $A \in
\ZZ^{d\times n}$ is given by \begin{equation} \vol_{\ZZ A} (A)=d!
\dfrac{\vol_{\RR^d}(\Delta_A )}{[\ZZ^d : \ZZ A]} \label{volume}
\end{equation} where $[\ZZ^d : \ZZ A]$ is the index of the subgroup $\Z A :=\sum_{i=1}^n
\Z a_i \subseteq \Z^d$, $\Delta_A$ is the convex hull of the columns
of $A$ and the origin in $\RR^d$ and $\vol_{\RR^d}(\Delta_A)$
denotes the Euclidean volume of the polytope $\Delta_A$.

Let us also recall that the direct sum of two matrices $A_1 \in \ZZ^{d_1 \times n_1}, A_2
\in \ZZ^{d_2 \times n_2}$ is the following $(d_1 + d_2 )\times (n_1
+ n_2)$ matrix:
$$
A_1 \oplus A_2 = \left(
\begin{array}{cc}
                                                        A_1 & 0_{d_1 \times n_2} \\
                                                        0_{d_2 \times n_1 } & A_2 \end{array}
\right) $$ where $0_{d\times n}$ denotes the $d\times n$ zero matrix.

The following two lemmas are easy to prove.

\begin{lemma}\label{volume-direct-sum}
If $A$ is the direct sum of two matrices $A_1 \in \ZZ^{d_1 \times
n_1}, A_2 \in \ZZ^{d_2 \times n_2}$ then $\vol_{\ZZ A} (A)=
\vol_{\ZZ A_1}(A_1) \cdot \vol_{\ZZ A_2}(A_2)$.
\end{lemma}

\begin{lemma}\label{lemma-disjoint-var} Let $A_i \in \ZZ^{d_i \times n_i}$ be full rank matrices, $d_i\leq n_i$,
and $\beta_{(i)}\in \CC^{d_i}$ for $i=1,2$. If $A=A_1 \oplus A_2$
and $\beta =(\beta_{(1)},\beta_{(2)})$ then we have that $H_A (\beta)=D H_{A_1}(\beta_{(1)} ) + D
H_{A_2}(\beta_{(2)} )$ where $H_{A_1}(\beta_{(1)})$ is a left ideal
of the Weyl Algebra $D_{A_1}=\CC [x_1 ,\ldots ,x_{n_1}]\langle \partial_{1} ,\ldots ,\partial_{n_1 } \rangle$ and
$H_{A_2}(\beta_{(2)})$ is a left ideal of the Weyl Algebra
$D_{A_2}=\CC [x_{n_1 +1} ,\ldots ,x_{n_1 +n_2}]\langle \partial_{n_1 +1} ,\ldots ,\partial_{n_1 +n_2} \rangle$ (equivalently, $M_A
(\beta )$ is the exterior tensor product of $M_{A_1}(\beta_{(1)} )$
and $M_{A_2} (\beta_{(2)})$).

\end{lemma}

The following corollary follows from Lemma \ref{lemma-disjoint-var}
by general properties of the exterior tensor product of holonomic
$D$--modules.

\begin{cor}\label{corollary-product-rank} Under the assumptions of Lemma \ref{lemma-disjoint-var}
we have:
\begin{enumerate}
\item[i)] $\rank (M_A (\beta ))= \rank (M_{A_1} (\beta_{(1)}))\cdot \rank
(M_{A_2} (\beta_{(2)})) $.
\item[ii)] If $\Omega_i$ is a basis for the space of (holomorphic) solutions of
the hypergeometric system $M_{A_i} (\beta_{(i)})$ at a point $p_i
\in \CC^{n_i}$, then the set $$\Omega = \{f_1 (x_1 ,\ldots ,
x_{n_1}) \cdot f_2 (x_{n_1 +1} ,\ldots ,x_{n_1 + n_2}) : f_i \in
\Omega_i , i=1,2\}$$ is a basis for the space of (holomorphic)
solutions of $M_{A}(\beta )$ at $p=(p_1 , p_2 )\in \CC^{n_1 + n_2}$.
\end{enumerate}
\end{cor}

In view of Corollary \ref{corollary-product-rank}, we can already give a first type of families of hypergeometric systems for which the rank jump grows exponentially with $d$.

\begin{theorem}\label{family-product}
Let $A\in \ZZ^{d\times n}$ and $\beta \in \CC^d$ be such that $M_A (\beta)$ has a rank jump, i. e. $\rank (M_{A} (\beta ))/ \vol_{\ZZ A}(A) = q >1$. 
Consider for $d_r=r d$ with $r\geq 1$ the matrix $A_{r}\in \ZZ^{d_r \times n_r } $ ($n_r =rn$), defined as the direct sum of $r$ copies of $A$, and the parameter vector 
$\beta_{r}=(\beta , \ldots ,\beta ) \in \CC^{d_r}$, defined by $r$ copies of $\beta$ as well. We have that the family given by 
$(A_{r},\beta_r )$ satisfies $\rank (M_{A_r} (\beta_r ))/\vol_{\ZZ A_r }(A_r ) \geq a^{d_r}$ where $a=\sqrt[d]{q}>1$.
\end{theorem}

In the sequel, we will first consider an example of a family similar to the ones given by Theorem \ref{family-product} and then, we will modify this example 
in order to exhibit for all $d\geq 2$ a family of $A$--hypergeometric systems with exponential growth of rank jumps on $d=\rank (A)$ which are not exterior tensor products of smaller hypergeometric systems.

\begin{example}\label{rank-jump-d2}
For $d=2$ we will consider the first example of a hypergeometric system with
rank jump described in \cite{ST98}. Consider the pair $(A_{(2)} ,\beta_{(2)})$
where
\begin{equation} A_{(2)}=\left(\begin{array}{cccc}
                                           1 & 1 & 1 &  1 \\
                                           0 & 1 & 3 & 4
                                         \end{array}\right)\label{first-example} \mbox{ and } \beta_{(2)} = \displaystyle\left(\begin{array}{c}
1\\
2\end{array}\right). \end{equation} The toric ideal associated with $A_{(2)}$ is
$$I_{A_{(2)}}=\left( \partial_1 \partial_4 -\partial_2 \partial_3 , \partial_1^2 \partial_3 -\partial_2^3 ,
\partial_2 \partial_4^2 -\partial_3^3 , \partial_1 \partial_3^2 -\partial_2^2 \partial_4\right)$$ and the Euler operators
are $E_1 -\beta_{(2),1} = x_1 \partial_1 + x_2 \partial_2 + x_3 \partial_3 +x_4 \partial_4 -1$ and
$E_2 -\beta_{(2),2} =  x_2 \partial_2 + 3 x_3 \partial_3 + 4 x_4 \partial_4 -2 $.

For this example $\rank ( M_{A_{(2)}}(\beta ))=\vol_{\ZZ A_{(2)}}(A_{(2)} )=4$ for all 
$\beta \in \CC^2 \setminus \{\beta_{(2)}\}$ but $\rank ( M_{A_{(2)}}(\beta_{(2)} ))=5$. A basis of the space of solutions of $M_{A_{(2)}}(\beta_{(2)} )$ 
can also be found in \cite{ST98}. Let us point out that this basis consists of the two Laurent polynomials $p_1  =x_2^2 /x_1$, 
$p_2 = x_3^2 / x_4$ and other $3$ functions that are not Laurent polynomials.
\end{example}

\begin{example}\label{rank-jump-d3}
 For $d=3$, we will consider the hypergeometric system of the family $\{ M_{A_{(d)}} (\beta_{(d)})\}_{d\geq 2}$
described in \cite{MW}. It is the one associated with the pair
\begin{equation}
A_{(3)}=\left(\begin{array}{cccccc}
 1 & 1 & 1 & 1 & 1 & 1 \\
 0 & 0 & 0 & 0 & 1 & 1 \\
 0 & 1 & 3 & 4 & 0 & 1 \end{array}\right) \mbox{ and } \beta_{(3)} =\left( \begin{array}{c}
                                                                                      1\\
0\\
2
                                                                                     \end{array}\right) \label{rank-jump-3}
\end{equation} The volume of $A_{(3)}$ is $d+2=5$ while the rank of $M_{A_{(3)}}(\beta (3))$ is $2d+1=7$.

\end{example}

\begin{example}\label{rank-jump-d-product}
For any $d\geq 4$, let $r,s \in \NN$ be such that $2r+3s=d$. We will choose $s$ as high as possible in order to 
fix uniques $r,s\in \NN$ for each $d$ (in particular $0\leq r \leq 4$). 

We define $A_{(d)} \in \ZZ^{d\times 2d}$ to be the direct sum of $r$ copies of the matrix $A_{(2)}$ and $s$ copies of the matrix $A_{(3)}$. By Lemma
\ref{volume-direct-sum} and examples \ref{rank-jump-d2} and \ref{rank-jump-d3} we have that $\vol_{\ZZ A_{(d)} } (A_{(d)})=4^r 5^s$.

On the other hand, let $\beta_{(d)} \in \CC^d$ be the complex vector with coordinates 
$\beta_{(d),2i-1}=1$ and $\beta_{(d),2i}=2$ for $1\leq i \leq r$ and
$\beta_{(d),2 r+ 3j-2}=1, \beta_{(d),2r+ 3 j-1}=0, \beta_{(d), 2r+3j}=2$ for $1
\leq j \leq s$ (i. e., $\beta_{(d)}$ has a copy of $\beta_{(2)}$ for each
copy of $A_{(2)}$ and a copy of $\beta_{(3)}$ for each copy of
$A_{(3)}$). With this definition of $(A_{(d)} ,\beta_{(d)} )$ and using
Corollary \ref{corollary-product-rank} and examples \ref{rank-jump-d2} and \ref{rank-jump-d3} we have that $\rank
( M_{A_{(d)}} (\beta_{(d)} )) = 5^r 7^s$. Thus 
$\rank ( M_{A_{(d)}}(\beta (d)))/\vol_{\ZZ A_{(d)}}(A_{(d)})=(5/4)^r (7/5)^s \geq (\sqrt{5}/2)^d$.
\end{example}

\begin{remark}\label{remark-polynomial-solutions}
Example \ref{rank-jump-d-product} also shows that the rank jump $j_A (\beta )$ can be greater than the number of 
Laurent polynomial solutions of $M_A (\beta)$. Indeed, since the space of Laurent polynomial solutions of $M_{A_{(2)}}(\beta_{(2)})$ has dimension $2$ 
(see \cite{ST98}) and the space of Laurent polynomial solutions of $M_{A_{(3)}}(\beta_{(3)})$ has dimension $4$ (see \cite{MW}) then, by Corollary \ref{corollary-product-rank}, the space of 
of Laurent polynomial solutions of $M_{A_{(d)}}(\beta (d))$ has dimension $2^r 4^s< j_{A_{(d)}}(\beta_{(d)})=5^r 7^s - 4^r 5^s$ for $r,s\geq1$.
\end{remark}

We are going to modify Example \ref{rank-jump-d-product} in order to get
hypergeometric systems that are not exterior tensor products of smaller hypergeometric systems.

Consider the following matrices and parameters:

\begin{equation} \hat{A}_{(2)}=\left(\begin{array}{ccccc}
                                           1 & 2 & 2 & 2 & 2 \\
                                           0 & 0 & 1 & 3 & 4
                                         \end{array}\right)\label{base1} \mbox{ and } \hat{\beta}_{(2)} = \displaystyle\left(\begin{array}{c}
3\\
2\end{array}\right). \end{equation}

 \begin{equation}
\hat{A}_{(3)}=\left(\begin{array}{ccccccc}
 1 & 2 & 2 & 2 & 2 & 2 & 2 \\
 0 & 0 & 0 & 0 & 0 & 1 & 1 \\
 0 & 0 & 1 & 3 & 4 & 0 & 1 \end{array}\right) \mbox{ and parameter } \hat{\beta}_{(3)} =\left( \begin{array}{c}
                                                                                      3\\
0\\
2
                                                                                     \end{array}\right) \label{base2}
\end{equation} 

Notice that $\hat{A}_{(2)}$ and $\hat{A}_{(2)}$ are obtained from $A_{(2)}$ and $A_{(3)}$ respectively by multiplying the first 
row by $2$ (this doesn't change the hypergeometric system) and then by adding a first column with its first coordinate equal to 
$1$ and the other coordinates equal to zero. After these modifications we get that $\vol_{\ZZ \hat{A}_{(2)}}(\hat{A}_{(2)})=2 \cdot 4 =8$ 
and that $\vol_{\ZZ \hat{A}_{(3)}}(\hat{A}_{(3)})=2 \cdot 5 =10$. However, since $\hat{\beta}_{(i)}$ is a 
\emph{hole} in $\NN \hat{A}_{(i)}$ (meaning that $\hat{\beta}_{(i)}\notin \NN \hat{A}_{(i)}$ but $\hat{\beta}_{(i)} +(\NN \hat{A}_{(i)}\setminus \{0\}) \subseteq \NN \hat{A}_{(i)}$) we have by Remark 4.14 in \cite{Oku} 
that $\rank (M_{\hat{A}_{(i)}}(\hat{\beta}_{(i)}))=\vol_{\ZZ \hat{A}_{(i)}}(\hat{A}_{(i)})+(i-1)$, $i=2,3$.

The following Lemma follows from the results in \cite{Berkesch}.

\begin{lemma}\label{equal-rank}
Let $A \in \ZZ^{d\times n}$ and $B\in \ZZ^{d\times m}$ be two matrices verifying that $\NN A = \NN B$ and $\Delta_A =\Delta_B$ then
$\rank (M_A (\beta ))=\rank (M_B (\beta ))$ for all $\beta \in \CC^d$.
\end{lemma}

For $d=2r+3s\geq 2$, $r ,s \in \NN$ (with $s$ as high as possible), let $\hat{\beta}_{(d)} \in \CC^d $ be the complex vector that is given
by $r$ copies of  $\hat{\beta}_{(2)}$ and $s$ copies of $\hat{\beta}_{(3)}$. The new matrix $\hat{A}_{(d)} \in \ZZ^{d\times (6 r +8 s-1)}$ is
constructed as follows.

Let $a_1 , a_2 ,\ldots , a_{5 r +
7s} \in \ZZ^d$ be the columns of the matrix $A_{r,s}
=\hat{A}_{(2)}\oplus
\stackrel{\underbrace{r}}{\cdots}\oplus\hat{A}_{(2)}\oplus
\hat{A}_{(3)} \stackrel{\underbrace{s}}{\cdots} \oplus \hat{A}_{(3)}
\in \ZZ^{d \times (5r+7s)}$.

We will construct a matrix $\hat{A}_{(d)}$ by adding $r+s-1$ column vectors to the matrix $A_{r,s}$. 
This vectors will belong to both $\Delta_{A_{r,s}}$ and $\NN A_{r,s}$. 
These conditions guarantee that $\vol_{\ZZ \hat{A}_{(d)}} (\hat{A}_{(d)})=\vol_{\ZZ A_{r,s}} (A_{r,s})=8^r 10^s$ and by  Lemma 
\ref{equal-rank}, we will also have that 
$\rank (M_{\hat{A}_{(d)}} (\beta ))=\rank (M_{A_{r,s}} (\beta))$ for all $\beta \in \CC^d$. In particular, for $\beta =\hat{\beta}_{(d)}$, 
we have $\rank (M_{\hat{A}_{(d)}} (\hat{\beta}_{(d)} ))=9^r 12^s$.

If $r\geq 2$ then for $1\leq i \leq r-1$ we define: $$a_{5r+7s+i}=a_1
+ a_{5i+1}=\dfrac{1}{2} a_2 + \dfrac{1}{2} a_{5i+2} \in \NN A_{r,s}\cap \Delta_{A_{r,s}}.$$ Notice that $(a_{5r+7s+i})_j$ equals $1$ for $j=1,2i+1$
and $0$ otherwise.

If $r, s\geq 1$ then for $1\leq i \leq s$ we define
$$a_{5r+7s+r-1+i}=a_1 + a_{5r+7i+1}= \dfrac{1}{2} a_2 + \dfrac{1}{2} a_{5r+7i+2}\in \NN A_{r,s}\cap \Delta_{A_{r,s}}.$$

If $r=0$ and $s\geq 1$ then for $1\leq i \leq s-1$ we define
$$a_{7s+i}=a_1 + a_{7i+1}=\dfrac{1}{2} a_2 + \dfrac{1}{2} a_{7i+2}\in \NN A_{r,s}\cap \Delta_{A_{r,s}}.$$

Let us define $\hat{A}_{d}=(a_1 \; a_2 \; \ldots \; a_{6r+8s-1})$
and recall that $\hat{\beta}_{(d)}\in \CC^d$ is given by $r$ copies of $\hat{\beta}_{(2)}$
and $s$ copies of $\hat{\beta}_{(3)}$. The hypergeometric system $M_{\hat{A}_{(d)}} (\hat{\beta}_{(d)} )$ is not 
an exterior tensor product of smaller hypergeometric systems and we have proved the following.

\begin{theorem} With the notations above we have
$$\dfrac{\rank (M_{\hat{A}_{(d)}} (\hat{\beta}_{(d)} ))}{\vol_{\ZZ \hat{A}_{(d)}} (\hat{A}_{(d)})}=(9/8)^r (12/10)^s \geq (\sqrt{9/8})^d$$
\end{theorem}

\begin{remark}
Notice that the toric ideal associated with $\hat{A}_{(d)}$ is not
homogeneous. However, by Theorem 7.3 in \cite{Berthesis}, if we
consider the associated homogeneous matrix $\hat{A}_{(d)}^h$ (that
is obtained by adding to the matrix $\hat{A}_{(d)}$ a first column
of zeroes and after that a first row of ones) and the parameter
$\hat{\beta}_{(d)}^{h}=(\beta_0 , \hat{\beta}_{(d)})$ with $\beta_0
\in \CC$ then the rank of $M_{\hat{A}_{(d)}^{h}}(\hat{\beta}_{(d)})$
equals the rank of $M_{\hat{A}_{(d)}}(\hat{\beta}_{(d)})$ if
$\beta_0 \in \CC$ is generic. This implies that for particular
$\beta_0$ (for example $\beta_0 =0$) the rank of $M_{\hat{A}_{(d)}^{h}}(\hat{\beta}_{(d)}^{h})$
will be greater than or equal to the rank of
$M_{\hat{A}_{(d)}}(\hat{\beta}_{(d)})$ by the upper semi--continuity
of the rank \cite{MMW}. Moreover,
$\vol_{\hat{A}_{(d)}^{h}}(\hat{A}_{(d)}^{h})=\vol_{\hat{A}_{(d)}}(\hat{A}_{(d)})$.
\end{remark}

\end{document}